\documentclass[11pt]{amsart}
\pdfoutput=1
\usepackage{latexsym, graphicx, epsfig, amsmath, amsfonts,amssymb}
\usepackage{multirow}
 \usepackage{diagbox}
\usepackage{color}
 \usepackage[percent]{overpic}
\usepackage{algorithm}
\usepackage{algorithmicx}
\usepackage{algpseudocode}
\floatname{algorithm}{Algorithm}
\usepackage{threeparttable}
\usepackage{booktabs}

\usepackage{epstopdf,mathtools,bbm}
\usepackage[left=1.25in,right=1.25in,top=1in]{geometry}

\def\mb{\mathbf}
\def\R{\mathbb{R}}
\def\N{\mathcal{N}}
\def\NN{\mathcal{NN}}

\newcounter{Rownumber}

\title{An acceleration strategy for randomize-then-optimize sampling via deep neural networks}

\author{Liang Yan}
\thanks{School of Mathematics, Southeast University, Nanjing, China. Email: yanliang@seu.edu.cn. L. Yan is supported by NSF of China (No.11771081), the science challenge project (No. TZ2018001) and  Zhishan Young Scholar Program of SEU}

\author{Tao Zhou}
\thanks{LSEC, Institute of Computational Mathematics, Academy of Mathematics and Systems
Science, Chinese Academy of Sciences, Beijing 100190, China. Email: tzhou@lsec.cc.ac.cn. T. Zhou is partially supported  by the National Key R$\&$D Program of China(No. 2020YFA0712000), the NSF of China (under grant numbers 11822111, 11688101and 11731006), the science challenge project (No. TZ2018001),  the Strategic Priority Research Program of Chinese Academy of Sciences (No.  XDA25000404) and youth innovation promotion association (CAS)}

\date{March 6, 2020}
\begin{document}

\maketitle

\begin{abstract}
Randomize-then-optimize (RTO)  is widely used for sampling from posterior distributions in Bayesian inverse problems. However, RTO may be computationally intensive for complexity problems due to repetitive evaluations of the expensive forward model and its gradient.  In this work, we present a novel strategy to substantially reduce the computation burden of RTO by using a goal-oriented  deep neural networks (DNN) surrogate approach. In particular, the training points for the DNN-surrogate are drawn from a local approximated posterior distribution, and it is shown that the resulting algorithm can provide a flexible and efficient sampling algorithm, which converges to the direct RTO  approach. We present a Bayesian inverse problem governed by a benchmark elliptic PDE to demonstrate the computational accuracy and efficiency of our new algorithm (i.e., DNN-RTO). It is shown that with our algorithm, one can significantly outperform the traditional RTO.
\end{abstract}

%\begin{keywords}
%Stein variational gradient descent, Bayesian inversion, surrogate modeling, deep neural network.
%\end{keywords}

\pagestyle{myheadings}
\thispagestyle{plain}
\markboth{LIANG YAN, AND TAO ZHOU}
{AN ACCELERATION STRATEGY FOR RTO SAMPLING VIA DNN}

\maketitle

%%%% Start %%%%%%
\section{Introduction}
The Bayesian approach provides a systematic framework for quantifying the uncertainty in the parameter estimation for inverse problems \cite{Kaipio+Somersalo2005,Stuart2010}.  In the Bayesian approach, the prior knowledge of the unknown parameters and the forward model are combined to yield a {\it posterior} probability distribution.  Then the unknown parameters can be characterized by their posterior distributions. The main task of the Bayesian approach is to draw samples from the posterior distributions, and then evaluate the associate statistic information, e.g., expectation, variance, etc. Since analytical formulas of the posterior are in general not available, many numerical sampling approaches such as Markov chain Monte Carlo (MCMC) methods \cite{Brooks2011} have been developed.

Using MCMC for sampling from the posterior is often computationally challenging. Firstly, each evaluation of the system output involves a forward model evaluation, and this is infeasible if the model is expensive to evaluate. Secondly, the geometry of the posterior distribution may admits complex features in the parametric space (such as local concentration).  To reduce the computational complexity, one possible way is to use the so-called surrogate approach: instead of using the true forward problems, one constructs a surrogate to the true model and then samples from the posterior distribution induced by the surrogate. In case the surrogate is computationally less expensive, one can speed up the MCMC algorithms dramatically.  Different surrogate approaches have been investigated in recent years, for example, projection-type reduced order models \cite{Arridge2006,Galbally+Fidkowski+Willcox+Ghattas2010,Lieberman+Willcox+Ghattas2010},  polynomial chaos (PC) based surrogates \cite{Marzouk+Najm2009,Marzouk+Najm+Rahn2007,Marzouk+Xiu2009,yan+guo2015,Yan+Zhang2017IP}, and Gaussian process regression \cite{kennedy2001,Rasmussen2003,stuart+teckentrup2016}, to name a few.  Recent works also include multi-fidelity surrogate approach \cite{Yan+Zhou19JCP, Yan+Zhou2019ADNN}  where surrogate models will be adaptively refined  during the MCMC sampling procedure. See also an application of the multi-fidelity surrogate approach to derivative-free methodologies, e.g. Ensemble Kalman inversion \cite{Yan+Zhou2019PCEKI}.
Although the surrogate approach can be very effective when exploring low-dimensional distributions, it can be very inefficient for complex, high-dimensional distributions with local concentration\cite{Brooks2011}:  successive states may exhibit high autocorrelation, due to the random walk nature of the movement. As a result, the effective sample size (ESS) tends to be quite low and the convergence to the true distribution is usually very slow. To address this challenge, several strategies that use the geometry information of the posterior (such as the gradient, Hessian, and higher order derivatives) have been exploited to accelerate the convergence of MCMC, see, e.g. \cite{Beskos2017Geometric,Hoffman2014NUT,Girolami2011Riemann,Lan2016emulation,Martin2012SNMC}.

In the present work, we propose a new approach that combines a deep neural networks (DNN) surrogate and an optimization-based sampling approach for large-scale PDE constrained Bayesian inverse problems. In particular, we focus on the randomize-then-optimize (RTO) approach \cite{Wang+Cui2019scalable, Bardsley2014SISC,Wang+Bardsley2017SISC} which uses repeated solutions of a randomly perturbed optimization problem to produce samples from a non-Gaussian distribution (which is used as a Metropolis independence proposal).  Compared to the classical Metroplis-Hastings (MH) random walk algorithm, RTO admits higher acceptance probability and lower sample auto-correlation even for high-dimensional problems \cite{Wang+Cui2019scalable}. The computational complexity of the our approach depends on the structure of the corresponding optimization problem, which requires the evaluation of the forward model and its gradient. To this end, we construct a DNN-based surrogate which makes the optimization problems rather efficient to solve. More specifically, to obtain an accurate and efficient DNN-surrogate, we choose the training points from a local approximated posterior distribution, and this makes the training procedure very efficient. We summarize the main features of our approach in the following:
\begin{itemize}
\item A new approach that combines the RTO and a DNN surrogate (DNN-RTO). Notice that the DNN approach is a powerful tool for approximating high dimensional problems \cite{Han+Jentzen+E2018PNAS,Raissi2019JCP,Schwab+Zech2019AA,Tripathy+Bilionis2018JCP,Zhu+Zabaras2018bayesian,GeorgeI,GeorgeII,YING}.
\item To train the DNN-surrogate, we choose the training points from a local approximated posterior distribution. Thus the training procedure can be very efficient.
\item We present numerical examples to demonstrate that our DNN-RTO approach can achieve several order magnitude speedup, yet guarantee the same accuracy as the traditional RTO.
\end{itemize}

The rest of the paper is organized as follows. In the next section, we give a brief introduction to the Bayesian inverse problems.  In Section 3, we introduce the RTO algorithm.  Details of our new approach are presented in Section 4. In Section 5, we use two nonlinear inverse problems to demonstrate the accuracy and efficiency of the new approach. Finally, we give some concluding remarks in Section 6.

\section {Bayesian inverse problems}\label{sec:setup} We are interested in the problem of estimating an unknown parameter $u\in \R^n$ from indirect observations $d\in \R^m$ via the following forward model
\begin{eqnarray}\label{feq}
d = F(u)+e,
\end{eqnarray}
where  $F: \R^n \rightarrow \R^m$ is a parameter-to-observation map that maps the unknown parameter $u$ to the measurements $d$, and $e \sim\N(0,\Gamma_{\text{obs}}) $ is the mean-zero Gaussian noise with a symmetric positive definite covariance $\Gamma_{\text{obs}} \in \R^{m\times m}$.

In the Bayesian setting, the prior belief about the parameter $u$ is encoded in the prior probability distribution $\pi(u)$. Moreover, we assume that $u$ is independent of the noise $e$. The aim of the Bayesian inverse problem is  to infer the distribution of $u$  conditioned on the data $d$, i.e., the {\it posterior} distribution $\pi(u|d)$.  By the Bayes' rule, we have
\begin{eqnarray}\label{ppdf}
\pi(u|d) \propto \exp(-\eta(u;d)) \pi(u),
\end{eqnarray}
where the term  $\exp(-\eta(u;d))$ is called {\it likelihood}, 
\begin{equation}\label{poteneq}
\eta(u;d) =\frac{1}{2} \|d-F(u)\|^2_{\Gamma_{\text{obs}}} = \frac{1}{2}(y-F(u))^T\Gamma_{\text{obs}}^{-1}(y-F(u))
\end{equation}
is a so-called {\it potential}, and $\propto$ denotes proportionality up to a scaling constant that depends on $d$ (but not on $u$).

The central task of Bayesian inverse problems (BIPs) is to characterize  the posterior distribution (\ref{ppdf}), e.g., computing certain posterior statistic moments.  Notice that if the forward model $F$  is  nonlinear, then in general the potential yields a posterior distribution which cannot be written in a closed form. Consequently, standard sampling methods such as the Metropolis-Hastings (MH) sampler, have been extensively studied to sampling with the posterior distribution. In this work, we will focus on the randomize-then-optimize (RTO)-MH approach\cite{Bardsley2014SISC}. RTO-MH uses repeated solutions of a randomly perturbed optimization problem to produce samples from a non-Gaussian distribution, and this is used as a Metropolis independence proposal. It is noticed that RTO-NM often yields better results than traditional MCMC. In the next section, we shall provide with more details on the RTO-MH algorithm.

\section{RTO-Metropolis-Hastings algorithm}\label{sec:method}
For BIPs, the original RTO algorithm can only be used to sample from the posterior distribution when the prior distributions are Gaussian \cite{ Bardsley2014SISC}. This limitation was relaxed by transforming non-Gaussian prior densities into Gaussian densities, see \cite{Wang+Bardsley2017SISC} for example. We refer to \cite{Wang+Cui2019scalable}  for excellent introductions to RTO.  In this work, we shall focus on the Gaussian prior case, and follow closely the derivation in \cite{Wang+Cui2019scalable}.

We assume that the prior of the parameter $u$ is Gaussian, i.e., $u\sim \N(u_{\text{pr}},\Gamma_{\text{pr}})$.  Here $u_{\text{pr}}$ is the prior mean, and $\Gamma_{\text{pr}}$ is the prior covariance matrix.  Combining with Eq. (\ref{poteneq}), the posterior can be thus written as
\begin{equation}
\pi(u|d) \propto \exp\Big(-\frac{1}{2}\big(\|d-F(u)\|^2_{\Gamma_{\text{obs}}}+\|u-u_{\text{pr}}\|^2_{\Gamma_{\text{pr}}}\big)\Big).
\end{equation}
%{\color{red}Since the RTO requires that the posterior distribution admits a specific form, we need whiten the prior and the error model firstly. {\color{blue}Delete?}}
Using matrix factorizations of the covariances of prior and observation noise
\begin{equation*}
S_{\text{pr}}S_{\text{pr}}^T : = \Gamma_{\text{pr}}, \quad S_{\text{obs}}S_{\text{obs}}^T: = \Gamma_{\text{obs}},
\end{equation*}
we define a new variable $v: = S^{-1}_{\text{pr}}(u-u_{\text{pr}})$ which satisfies
\begin{equation}
0 = f(v)+\epsilon, \quad \epsilon \sim \N(0, I_m), \quad v \sim \N(0, I_n).
\end{equation}
Here $f(v)=S^{-1}_{\text{obs}}[F(S_{\text{pr}}v+u_{\text{pr}})-d]$,  $I_n$ and $I_m$ are identity matrices of size $n$ and $m$, respectively. Then, the resulting posterior density of the $v$ is given by
\begin{equation}\label{tarpdf}
\pi(v|d) = \pi_{\text{tar}}(v) \propto \exp(-\frac{1}{2}\|H(v)\|^2),
\end{equation}
where $H: \R^{n} \rightarrow \R^{(n+m)}$ is defined as
\begin{eqnarray}\label{eg1case1}
H(v) =\begin{bmatrix}
v \\
f(v)
\end{bmatrix}.
\end{eqnarray}

Notice that if we have a sample $v$ from the target density $\pi_{\text{tar}}(v)$, we can easily obtain the corresponding posterior sample of $u$ by applying the transformation $u=S_{\text{pr}}v+u_{\text{pr}}$. Now we outline how to use  the RTO-MH to sample from a posterior of the form (\ref{tarpdf}):
\begin{itemize}
\item We first choose a linearization point $v_{\text{ref}}$ from the following optimization problem
\begin{equation}\label{vrefpoint}
v_{\text{ref}} = \arg \min _{v} \frac{1}{2}\|H(v)\|^2.
\end{equation}
Then we can compute a matrix $Q\in \R^{(n+m)\times n}$ with orthonormal columns from a thin QR factorization of $\nabla H(v_{\text{\text{ref}}})$.
\item Next, we draw independent samples $\xi^{(i)}$ from an $n$-dimensional standard Gaussian, and for each sample $\xi^{(i)}$ we generate proposal points  $v^{(i)}_{\text{prop}}$  by solving the following optimization problem
\begin{equation}\label{orgopt}
v^{(i)}_{\text{prop}} = \arg\min_{v} \frac{1}{2}\|\bar{Q}^TH(v)-\xi^{(i)}\|^2,
\end{equation}
  The above equation (\ref{orgopt}) is called randomize-then-optimize (RTO), and under certain conditions (such as those in \cite{Bardsley2014SISC}), the points $v^{(i)}_{\text{prop}}$ are distributed according to the following proposal density,
\begin{equation}\label{orgRTOpro}
\pi_{\text{RTO}}(v) = (2\pi)^{-\frac{\pi}{2}}|\det(Q^T \nabla H(v))|\exp\Big(-\frac{1}{2}\|Q^TH(v)\|^2\Big).
\end{equation}
\item Finally, we use (\ref{orgopt}) and (\ref{orgRTOpro}) as an independence proposal within the MH algorithm for sampling from $\pi_{\text{tar}}(v)$. Given a previous sample $v^{(i-1)}$ and a proposed RTO sample  $v^{(i)}_{\text{prop}}$, the acceptance ratio for the MH method is given by \cite{Bardsley2014SISC}
\begin{equation}
\frac{\pi_{\text{tar}}(v^{(i)}_{\text{prop}})\pi(v^{(i-1)})}{\pi_{\text{tar}}(v^{(i-1)})\pi(v^{(i)}_{\text{prop}})} = \frac{w(v^{(i)}_{\text{prop}})}{w(v^{(i-1)})},
\end{equation}
where
\begin{equation}\label{wfun}
w(v) = |\det(Q^T\nabla H(v))|^{-1} \exp\Big(-\frac{1}{2}\|H(v)\|^2+\frac{1}{2}\|Q^TH(v)\|^2\Big).
\end{equation}
\end{itemize}
When the dimension of the parameter vector $v$ is very high, solving the optimization problem (\ref{orgopt}) and computing the RTO probability density (\ref{orgRTOpro}) can be computationally costly.  In order to overcome this challenge, Bardsley et.al. \cite{Wang+Cui2019scalable}  introduce a new subspace acceleration strategy to make the computational complexity of RTO scale linearly with the parametric dimension. The main idea is to use the singular value decomposition (SVD) of the linearized forward model $\nabla f(v_\text{ref})$ instead of computing the QR of the matrix $\nabla H(v_{\text{ref}})$. Similar to the original RTO, the scalable implementation of RTO also includes three steps:

\textbf{Step 1:}  Compute the  reduced SVD of the $\nabla f(v_\text{ref})$, which has rank $r$, as
\begin{equation}
\nabla f(v_\text{ref}) = \Psi \Lambda \Phi^T.
\end{equation}

\textbf{Step 2:}  Define
\begin{equation}
v_r = \Phi^Tv, \quad \mbox{and} \quad v = \Phi v_r+v_{\bot},
\end{equation}
where $v_{\bot}$ is an element in the orthogonal complement of range($\Phi$).  For each realization of an $n$-dimensional standard Gaussian random vector $\xi$, one can compute
\begin{equation}\label{vperps}
v_{\bot} = (I_n-\Phi\Phi^T)\xi
\end{equation}
and solve the following $r$-dimensional optimization problem with $\xi$ and $v_{\bot} $
\begin{equation}\label{scaopt}
v_r = \arg\min_{z} \|(\Lambda^2+I_r)^{-\frac{1}{2}}\big(z+\Lambda\Psi^T f(v_{\bot}+\Phi z)-\Phi^T \xi\big)\|^2.
\end{equation}
Notice that Eqs. (\ref{vperps}) and (\ref{scaopt}) are replaced the original $n$-dimensional optimization problem in (\ref{orgopt}).

\textbf{Step 3:} Compute the weighting function $w(v)$ in (\ref{wfun})  as
\begin{equation}\label{awfun}
w(v) = |\det(\tilde{Q}^T\nabla H(v))|^{-1} \exp\Big(-\frac{1}{2}\big(\|f(v)\|^2+\|\Phi^Tv\|^2-\|(\Lambda ^2+I_r)^{-\frac{1}{2}}(\Phi^Tv+\Lambda\Psi^Tf(v))\|^2\big)\Big),
\end{equation}
where the determinant takes the following simplified form
\begin{equation}\label{adetfun}
|\det(\tilde{Q}^T\nabla H(v))|=|\det(\Lambda^2+I_r)^{-\frac{1}{2}}||\det(I_r+\Lambda\Psi^T\nabla f(v)\Phi)|.
\end{equation}

The resulting MCMC method, which is called scalable implementation of RTO-MH, is summarized in Algorithm \ref{alg:SRTO}. %Although this scalable implementation  reduces the size of the optimization problem and reduce the size of the matrix determinant calculation for each proposal,

Notice that the main bottleneck for using the RTO-MH approach is the repetitive evaluations of forward model and its gradient in the optimization procedure (\ref{scaopt}) and the weighting function (\ref{awfun}) that involves the observation that may involve a complicated model. Thus it is desirable to construct effective approximation of these quantities to provides a good balance between accuracy and computation cost.  The motivates our DNN-based surrogate model approach which will be introduced in Section \ref{sec:dnnrto}.

\begin{algorithm}[t]
  \caption{Scalable implementation of RTO-MH \cite{Wang+Cui2019scalable}}
  \label{alg:SRTO}
  \begin{algorithmic}[1]
\State Find $v_{\text{ref}}$ using (\ref{vrefpoint})  .
 \State  Determine the Jacobian matrix of the forward model, $\nabla f(v_{\text{ref}})$.
 \State Compute the SVD of $\nabla f(v_{\text{ref}})$
    \For {$j=1,\cdots, n_{\text{samps}}$} {in parallel}
    \State  Samples $\xi^{(i)}$ from a standard $n$-dimensional Gaussian distribution.
    \State Solve for a proposal sample $v^{(i)}_{\text{prop}} = v_{\bot}+\Phi v_r$ using (\ref{vperps}) and (\ref{scaopt}).
    \State Compute $w(v^{(i)}_{\text{prop}})$ from (\ref{awfun}) using the determinant  from (\ref{adetfun}).
    \EndFor
 \State Set $v^0= v_{\text{ref}}$
 \For   {$j=1,\cdots, n_{\text{samps}}$} {in series}
 \State Sample $t$ from a uniform distribution on $[0,1]$.
      \State Accept/reject $v^{(i)}_{\text{prop}}$ according   $t< w(v^{(i)}_{\text{prop}})/w(v^{(i-1)}) $
 \EndFor
%     \EndProcedure
  \end{algorithmic}
\end{algorithm}
\section{DNN surrogate for RTO-MH}\label{sec:dnnrto}

In this section, we shall present a  DNN-based surrogate modeling to accelerate the RTO-MH approach.

\subsection{Feedforward DNN-based surrogate modeling}\label{sec:DNN}

\begin{figure}
\begin{center}
  \begin{overpic}[width=.65\textwidth,trim=20 0 20 15, clip=true,tics=10]{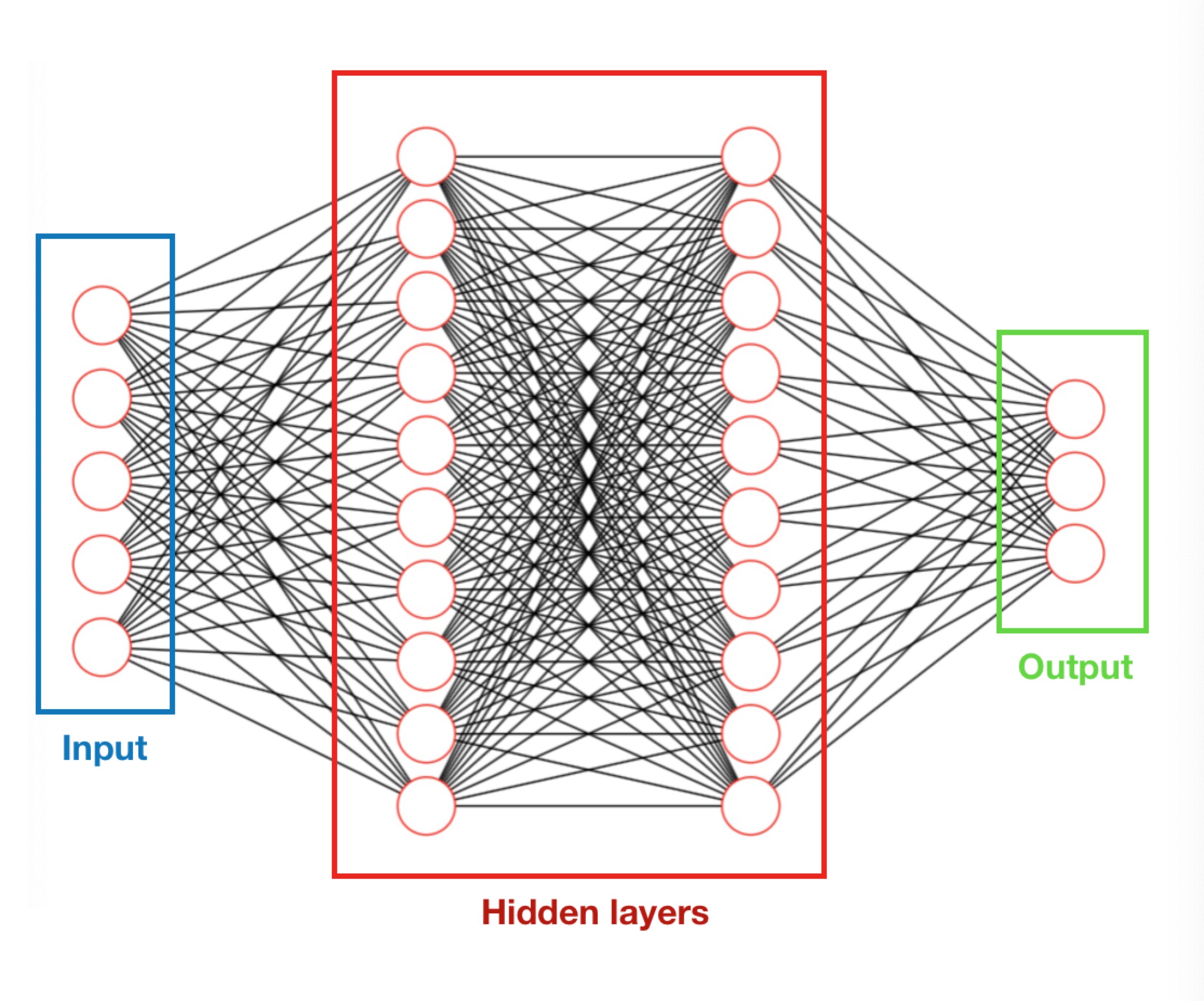}
  \end{overpic}
  \end{center}
\caption{The structure of a two-hidden-layer neural network.}\label{DNN_structure}
\end{figure}

The basic idea of deep neural networks (DNN) for surrogate modeling is that one can  approximate  an input-output map $f: \R^n\rightarrow \R^m$ through a hierarchical abstract layers of latent variables \cite{Goodfellow2016DL}. A typical example is the feedforward neural network, which is also called multi-layer perception (MLP). It consists of a collection of layers that include an input layer, an output layer, and a number of hidden layers.  Specifically, in the $k$th hidden layer, $d_k$ number of neurons are included. Each hidden layer of the network receives an output $z^{(k-1)}\in \R^{d_{k-1}}$ from the previous layer where an affine transformation of the form
\begin{equation}
\mathcal{F}_k(z^{(k-1)})  = \mb{W}^{(k)}z^{(k-1)}+\mb{b}^{(k)},
\end{equation}
is performed. Here  $\mb{W}^{(k)} \in \R^{d_{k}\times d_{k-1}}, \, \mb{b}^{(k)} \in \R^{d_{k}}$ are the weights and biases of the $k$th layer. The nonlinear activation function $\sigma$ is applied to each component of the transformed vector before sending it as an input to the next layer. The activation function is an identity function after an output layer. Thus, the final neural network representation is given by the composition
\begin{equation}
\NN(v) = (\mathcal{F}_L\circ\sigma\circ\mathcal{F}_{L-1}\circ \cdots \circ\sigma\circ\mathcal{F}_1)(v),
\end{equation}
where the operator $\circ$ is the composition operator, $v=z^{(0)}$ is the input. A typical neural network architecture can be founded in Fig. \ref{DNN_structure}. 

Some popular choices for the activation function include \texttt{sigmoid,  hyperbolic tangent,  rectified linear unit (ReLU)},  to name a few \cite{Goodfellow2016DL,Ramachandran2017}. In the current work, we shall use \texttt{Swish} as the activation function \cite{Ramachandran2017,Tripathy+Bilionis2018JCP}:
\begin{equation*}
\sigma(z)=\frac{z}{1+\exp(-z)}.
\end{equation*}
Once the network architecture is defined, one can resort to optimization tools to find the unknown parameters $\theta =\{ \mb{W}^{(k)},  \mb{b}^{(k)}\}$ based on the training data. Precisely, let $\mathcal{D}: =\{(v_i, y_i)\}^N_{i=1}$ be a set of training data, we can define the following minimization problem:
\begin{equation}\label{thetastar}
\arg\min_{\theta} \frac{1}{N}\sum^N_{i=1} \|y_i-\mathcal{NN}(v_i;\theta)\|^2,
\end{equation}
where $\mathcal{J}(\theta; \mathcal{D}) =  \frac{1}{N}\sum^N_{i=1} \|y_i-\mathcal{NN}(v_i;\theta)\|^2$ is the so called loss function.  Solving this problem is generally achieved by the stochastic gradient descent (SGD) algorithm \cite{Bottou2010}.   SGD simply minimizes the function by taking a negative step along an estimate  of the gradient  $\nabla_{\theta} \mathcal{J}(\theta; \mathcal{D})$ at iteration $k$. The gradients are usually computed through backpropagation.  At each iteration, SGD updates the solution by
\begin{equation*}
\theta_{k+1}=\theta_k-\lambda \nabla_{\theta} \mathcal{J}(\theta; \mathcal{D}),
\end{equation*}
where $\lambda$ is the learning rate. Recent algorithms that offer adaptive learning rates are available, such as \texttt{Ada-Grad} \cite{Zeiler2012Adadelta}, \texttt{RMSProp} \cite{Tieleman+Hinton2012lecture} and \texttt{Adam} \cite{Kingma2014Adam}, ect.  The present work adopts \texttt{Adam} optimization algorithm, and we shall construct a DNN as a surrogate model for the true forward model in the BIPs.

\subsection{Choosing effective training points}

Notice that the accuracy of the neural network $\NN(v;\theta)$ approximation are clearly influenced by the choice of the training points. A naive way to choose the training data is to generate enough data over the whole prior distribution, however, this may loss the gain computational efficiency \cite{Yan+Zhou19JCP,Yan+Zhou2019ADNN}. Notice that our concern in BIPs is only the posterior distribution which may be located in a small region. Consequently, we only need to make sure that $\NN(v;\theta)$ is accurate enough in the posterior density region where no need to ensure its accuracy elsewhere. Nevertheless, the high-probability density region of the posterior hard to be identified until data are available.  In the following, we propose a goal-oriented technique based on the observation data and the gradient information to choose the local training points.

Consider the first-order Taylor approximation of the forward model $f(v)$ at $v_{\text{ref}}$, we have
\begin{equation}
f(v)  \approx f(v_{\text{ref}}) +\nabla f(v_{\text{ref}})(v-v_{\text{ref}}): = Av-b,
\end{equation}
where $A:=\nabla f(v_{\text{ref}}) $ and $b=Av_{\text{ref}}-f(v_{\text{ref}})$. Using this linearization approximation, we can define the {\it approximation} posterior of $\pi_{\text{tar}}(v)$ as
\begin{equation}
\pi_{\text{tar}}(v) \approx \widetilde{\pi}_{\text{tar}}(v)\propto
\exp(-\frac{1}{2}\widetilde{H}(v)),
\end{equation}
where $\widetilde{H}: \R^{n} \rightarrow \R^{(n+m)}$ is defined as
\begin{eqnarray}
\widetilde{H}(v) =\begin{bmatrix}
v \\
Av-b
\end{bmatrix}.
\end{eqnarray}
It is easy to verify that $\widetilde{\pi}_{\text{tar}} \sim \N(\widetilde{\mu}_{\text{pos}}, \widetilde{\Gamma}_{\text{pos}})$ with
\begin{equation}
\widetilde{\Gamma}_{\text{pos}} = (A^TA+I)^{-1}, \quad \mbox{and} \quad \widetilde{\mu}_{\text{pos}}= \widetilde{\Gamma}_{\text{pos}} A^T b.
\end{equation}
This means that the first-order Taylor approximation provides us with a {\it local} Gaussian measure approximating $\pi_{\text{var}}$. This Gaussian measure allows for direct sampling if the reference point $v_{\text{ref}}$ and the gradient of the forward model are given.  Notice that the covariance is the inverse matrix of the $A^TA+I$, however, with a simple algebraic manipulation by the reduced SVD $A=\nabla f(v_{\text{ref}}) = \Psi \Lambda \Phi^T$, we can write it as
\begin{equation}\label{appcovmat}
\widetilde{\Gamma}_{\text{pos}} = \Phi(\Lambda^2+I_r)^{-1}\Phi^T+(I_n-\Phi\Phi^T).
\end{equation}
The above discussion motivates us to choose the training points from the local Gaussian approximation of the posterior with covariance $\widetilde{\Gamma}_{\text{pos}}$ defined in (\ref{appcovmat}). In Section \ref{sec:tests}, we will perform the comparison between this strategy and the prior-based strategy (choose the data in the whole prior domain) by numerical examples.

\subsection{DNN-based RTO-MH}
As mentioned in the previous sections, we can alleviate the complexity issue with a DNN surrogate for the forward model $f(v)$.  It is clear that after obtaining the parameters $\theta$, we have an explicit functional $\mathcal{NN}(v;\theta)$ and can compute its gradient $\nabla_v\mathcal{NN}(v;\theta)$ easily via the back propagation \cite{Bishop2006pattern}. These approximations can be then substituted into the computation procedure of the RTO-MH framework, and obtain the DNN-based RTO-MH algorithm. When the original forward model is computationally costly, simulating the surrogate $\NN(v;\theta)$ provides a more efficient mechanism for the RTO-MH sampler.  The detail of the scheme is summarized in Algorithm \ref{alg:DNN-RTO-MH}. Our proposed method provides a natural framework to incorporate DNN surrogate in RTO-MH. Moreover, it can be easily extended to other optimization-based sampling, e.g., the random-map implicit sampling \cite{Morzfeld2012random}.

\begin{algorithm}[th]
  \caption{The offline and online stages for the NN-RTO method}
  \label{alg:DNN-RTO-MH}
  \begin{algorithmic}[1]
  \State  {\bf Offline stage:}
  \State Choose $N$ training points $\{v_i\}^N_{i=1}$ randomly from the approximate posterior  $\widetilde{\pi}_{\text{var}}(v)$
   \State Compute the corresponding  full-order snapshots $\{f(v_i)\}^N_{i=1}$ ;
   \State Prepare the training set $\mathcal{D} = \Big\{v_i, f(v_i)\Big\}$;
    \State Train the DNN model $\NN(v;\theta)$ by using the training set $\mathcal{D}$.
    \vspace{0.3cm}
    \State {\bf Online stage:}
      \For {$j=1,\cdots, n_{\text{samps}}$} {in parallel}
    \State  Samples $\xi^{(i)}$ from a standard $n$-dimensional Gaussian distribution.
    \State Compute a proposal sample $v^{(i)}_{\text{prop}}$, and  $w(v^{(i)}_{\text{prop}})$ using RTO algorithm with the trained neural network $\NN(v;\theta)$.
    \EndFor
 \State Set $v^0= v_{\text{ref}}$
 \For   {$j=1,\cdots, n_{\text{samps}}$} {in series}
 \State Sample $t$ from a uniform distribution on $[0,1]$.
      \State Accept/reject $v^{(i)}_{\text{prop}}$ according   $t< w(v^{(i)}_{\text{prop}})/w(v^{(i-1)}) $
 \EndFor

%     \EndProcedure
  \end{algorithmic}
\end{algorithm}

\section{Numerical Examples}\label{sec:tests}

In this section, we present a benchmark elliptic PDE inverse problem to illustrate the accuracy and efficiency of the DNN-RTO-MH approach. To better present the results, we shall perform the following three-types of approaches:
\begin{itemize}
\item The conventional "RTO" (or the direct RTO) which is based on the true forward model evaluations.
\item The "NN-RTO-pr" that use a DNN surrogate with training data that are generated with respect to the prior distribution.
\item The "NN-RTO" that use a DNN surrogate with training data that are generated with respect to an approximation posterior distribution in Section 4.2, i.e., the suggested algorithm in this work.
\end{itemize}

In our all numerical tests, the computations were performed using MATLAB 2018a on an Intel-i7 desktop computer. For solving the optimization problem (\ref{scaopt}) we use the built-in nonlinear least squares solver \texttt{lsqnonlin} of MATLAB, which implements the trust region reflective Newton method.  The \texttt{Adam} optimizer is used to train the DNN as mentioned before. The learning rate is set to be $\lambda =5\times10^{-4}$, and the hyper-parameter values of \texttt{Adam} are chosen based on default recommendations as suggested in \cite{Kingma2014Adam}.
In order to compare our proposed method to standard RTO-MH in terms of sampling efficiency, we consider the effective sample size (ESS) adjusted by CPU time.  Given $n_{\text{samps}}$  posterior samples, the ESS for each parameters is defined as
\begin{equation}
ESS = \frac{n_{\text{samps}}}{1+2\sum^{K}_{k=1}\rho(k)}
\end{equation}
where $\sum^{K}_{k=1}\rho(k)$ is the sum of $K$ monotone sample autocorrelations\cite{Geyer1992practical}.
We use the minimum ESS over all parameters normalized by the CPU time, $s$(in seconds), as the overall measure of efficiency: min(ESS)/$s$.  In this work, we will use $n_{\text{samps}}=5000$ and $K=n_{\text{samps}}-1$ to obtain the results.

\subsection{Problem setup}

We consider the problem of inferring subsurface permeability from a finite number of noisy pressure head measurements \cite{Cui2014data,Yan+Zhou2019ADNN}.   The forward model is given by the solution of an elliptic PDE in two spatial dimensions
\begin{eqnarray}\label{2dellip}
-\nabla\cdot(\kappa(x) \nabla p(x))&=f(x),
\end{eqnarray}
where $x=(x_1,x_2)\in [0,1]^2$ is the spatial coordinate.
The boundary conditions are
\begin{align*}
& p(x)|_{x_1=0}=1, \quad p(x)|_{x_1=1}=0,\\
& \frac{\partial p(x)}{\partial x_2}|_{x_2=0}=x_1, \quad \frac{\partial p(x)}{\partial x_2}|_{x_2=1}=1-x_1.
\end{align*}
The data $d$ is given by a finite set of $p$, perturbed by noise, and the problem is to recover the permeability $\kappa(x)$ from these measurements.  In what follows, we choose the source $f(x) = 100\sin(\pi x_1)\sin(\pi x_2)$, and solve the equation (\ref{2dellip}) using a standard Galerkin finite element method with bilinear basis function on a uniform 40-by-40 grid.

\subsection{Example 1: a nine-dimensional inverse problem}

 \begin{figure}
\begin{center}
    \begin{overpic}[width=0.4\textwidth,trim= 20 0 0 25, clip=true,tics=10]{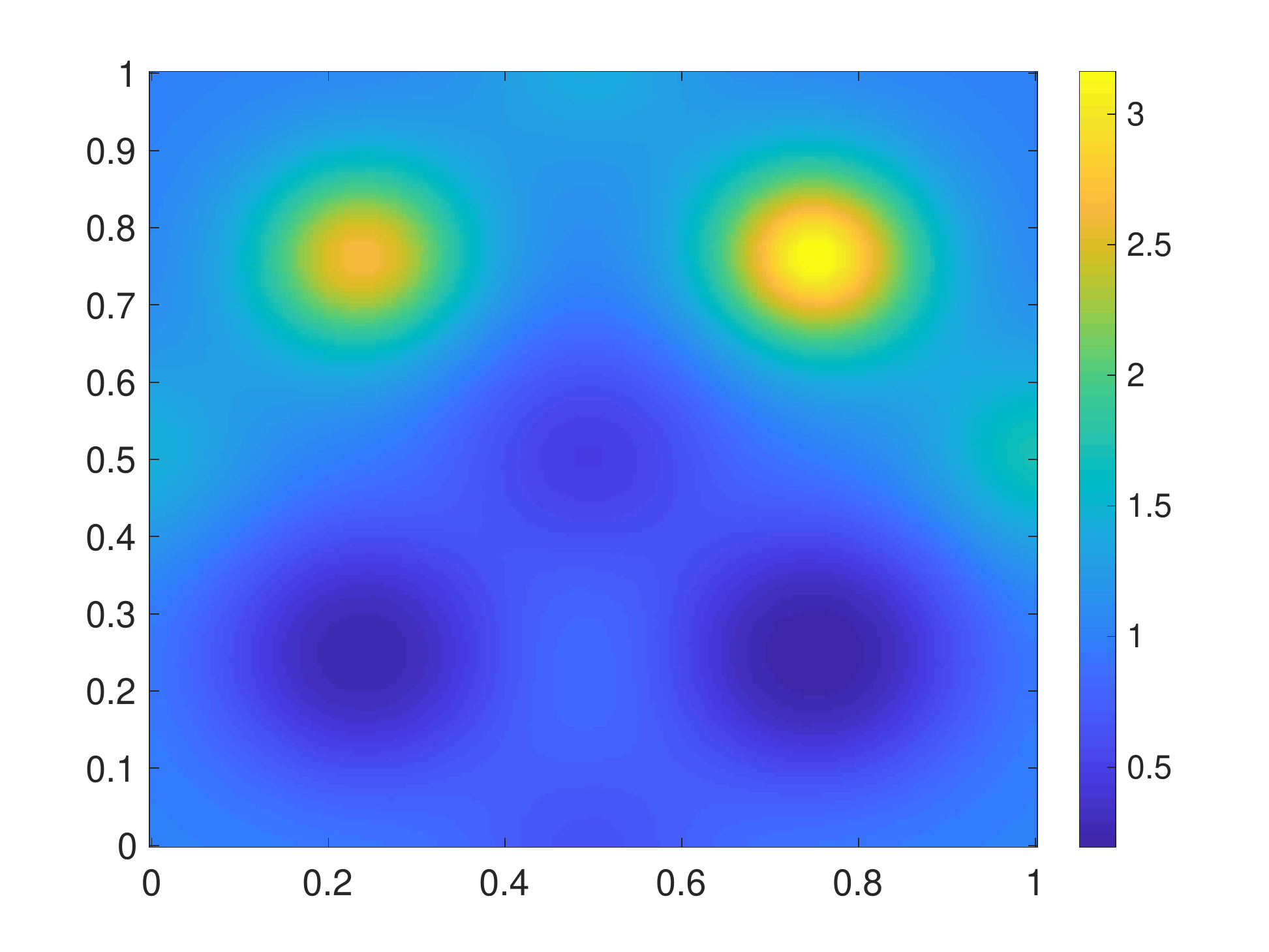}
   \end{overpic}
  \begin{overpic}[width=0.4\textwidth,trim= 20 0 0 25, clip=true,tics=10]{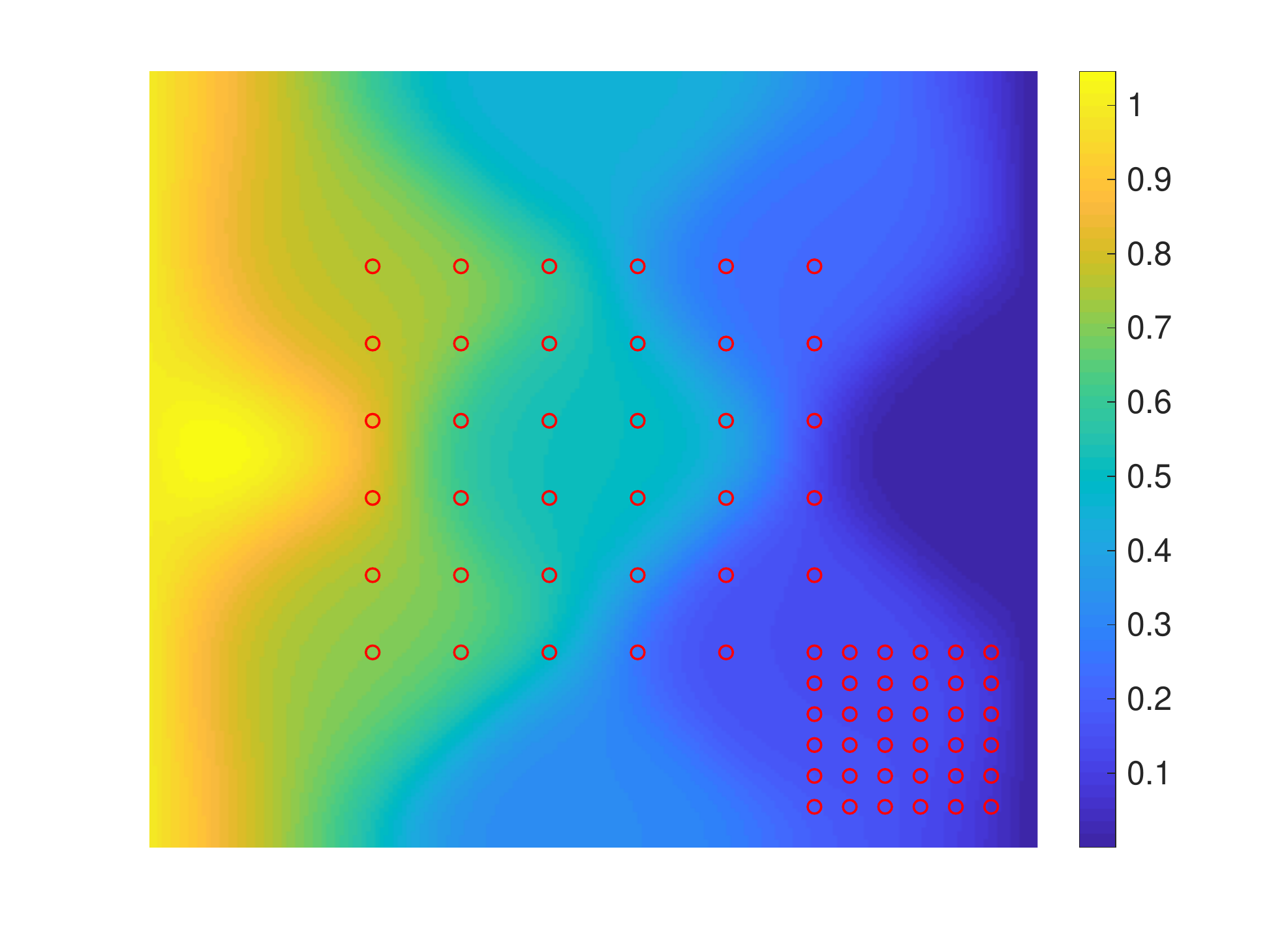}
 \end{overpic}
\end{center}
\caption{Example 1:  Setup of the test case for Example 1. Left: the true permeability used for generating the synthetic data sets. Right: the model outputs of the true permeability. The red circles indicate the measurement sensors.}\label{exact_eg1}
\end{figure}

In the first example, the permeability field $\kappa(x)$ is defined by
\begin{eqnarray*}
\kappa(x)=\sum^{9}_{i=1}\kappa_i \exp\left(- \frac{\|x-x_{0,i}\|^2}{2 \times 0.1^2}\right),
\end{eqnarray*}
where $\{x_{0,i}\}^{9}_{i=1}$ are the centers of the radial basis function, and the weights $\{\kappa_i\}^9_{i=1}$ are parameters in the Bayesian inverse problem.

This example is a typical benchmark problem considered in Refs. \cite{Cui2014data,Yan+Zhou19JCP}. We first choose a realization of $\kappa_i$ from a uniform distribution as the true solution.  The true permeability  field used to generate the test data and the corresponding pressure head are  shown in  Fig.\ref{exact_eg1}. The prior distributions on each of the weights $\kappa_i, i=1,\cdots, 9$ are independent and log-normal; that is, $\log(\kappa_i):=v^i\sim N(0,1)$.  Partial observations of the pressure field are collected 71 measurement sensors as shown by the red circles in Fig. \ref{exact_eg1}. This yields observed data $d\in \R^{71}$ as
\begin{eqnarray*}
d_j=p(x_j)+\max_{j}\{|p(x_j)|\}\delta\xi_j,
\end{eqnarray*}
where $\delta$ dictates the relative noise level and $\xi_j$  is a  Gaussian random variable with zero mean and unit standard deviation.  In the following, we set $\delta =0.05$.

\subsubsection{Computational efficiency}

 \begin{figure}
\begin{center}
    \begin{overpic}[width=0.7\textwidth,trim= 10 0 40 15, clip=true,tics=10]{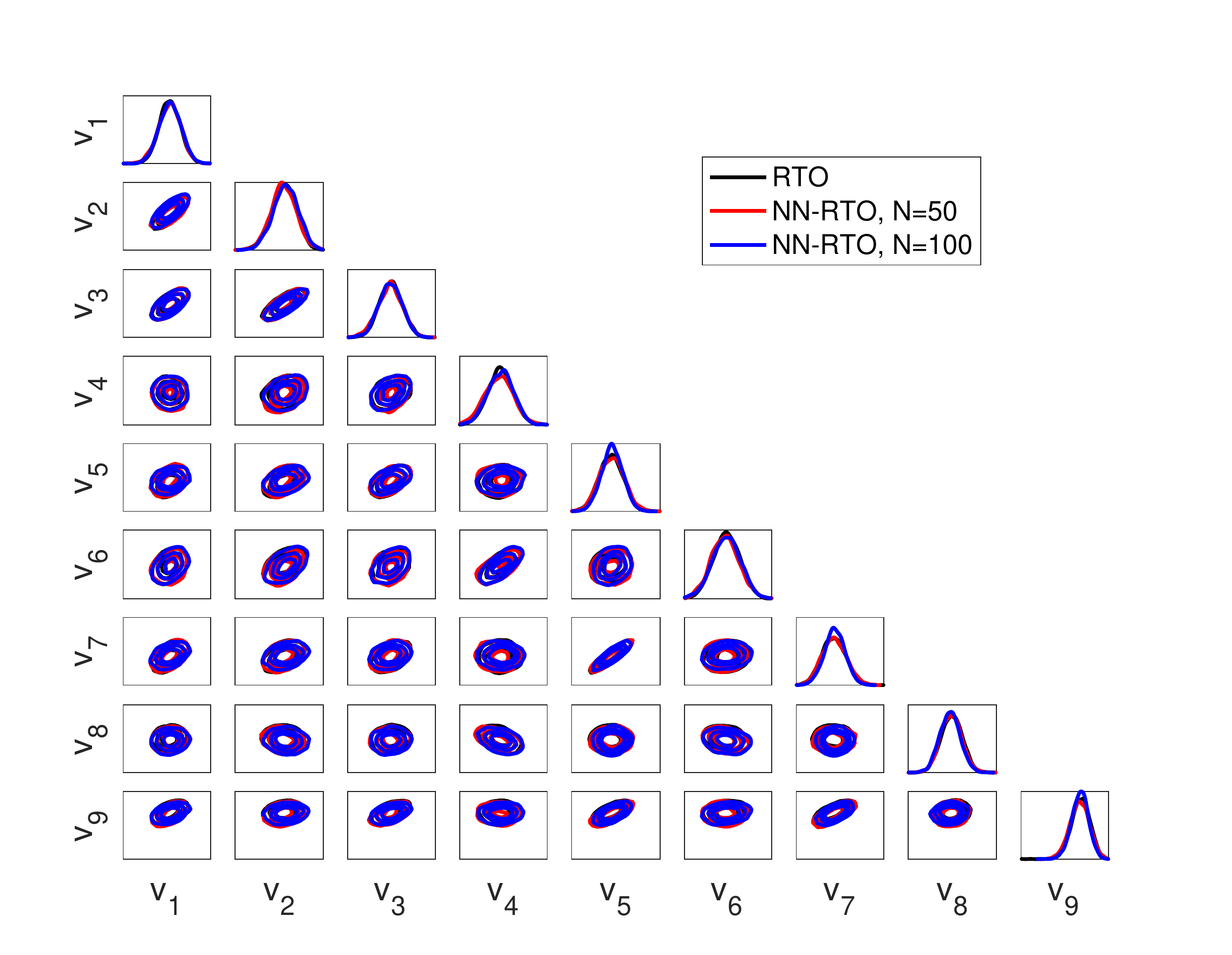}
    % \put (40,-3) {\scriptsize NN-RTO}
      \end{overpic}
%        \begin{overpic}[width=0.47\textwidth,trim= 10 0 40 15, clip=true,tics=10]{figures/resNNpr_setup0225_post-eps-converted-to.pdf}
%         \put (40,-3) {\scriptsize NN-RTO-pr}
%      \end{overpic}
\end{center}
\caption{One- and two-dimensional posterior marginals of the two parameters for Example 1 using NN-RTO approach. }\label{pos_contour_eg1}
\end{figure}

\begin{figure}
\begin{center}
    \begin{overpic}[height=6.4cm,width=4.25cm,trim= 35 10 45 15, clip=true,tics=10]{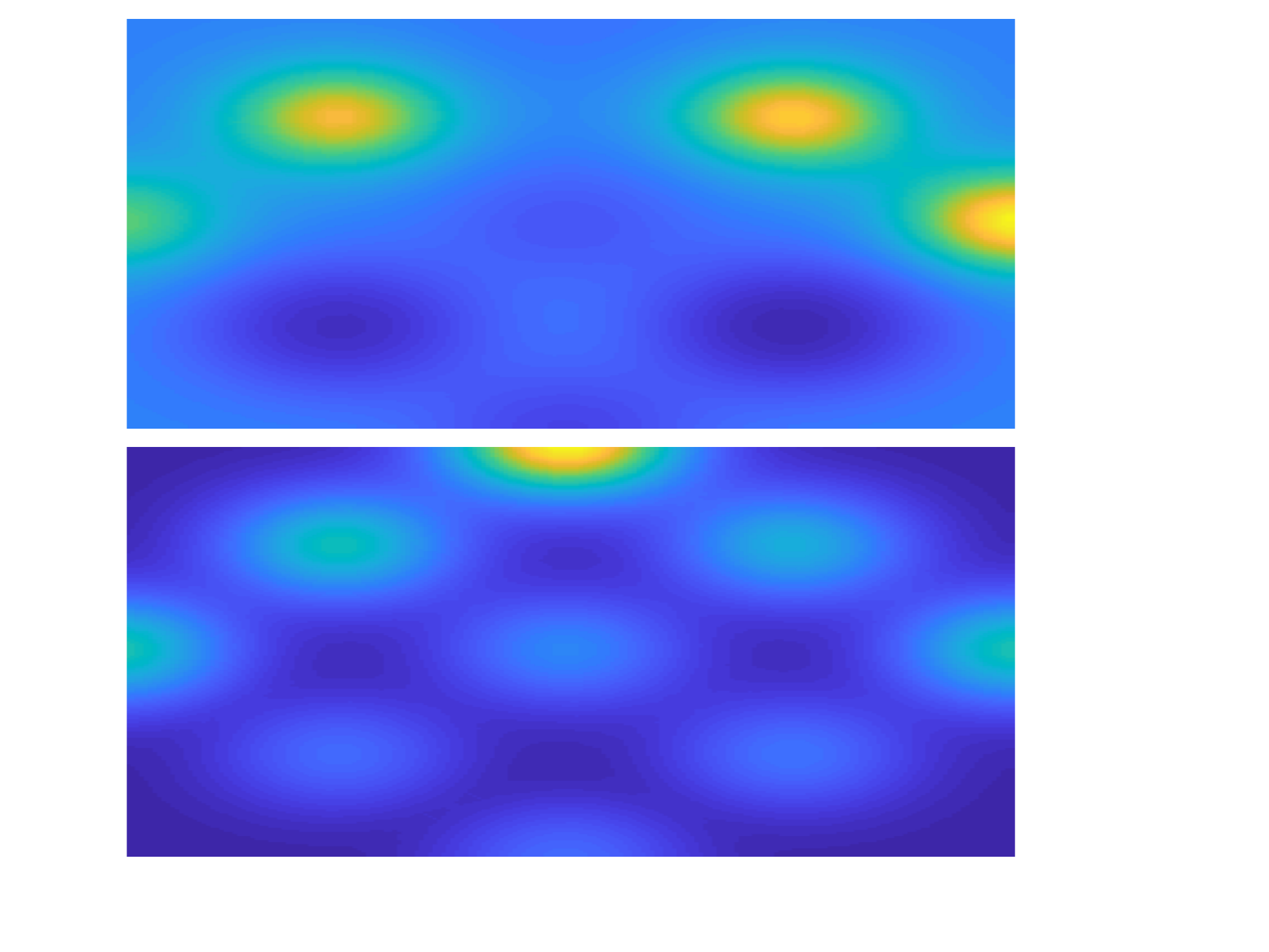}
         \put (23,-3) {\scriptsize RTO}
  \end{overpic}
    \begin{overpic}[height=6.4cm,width=4.25cm,trim= 35 10 45 15, clip=true,tics=10]{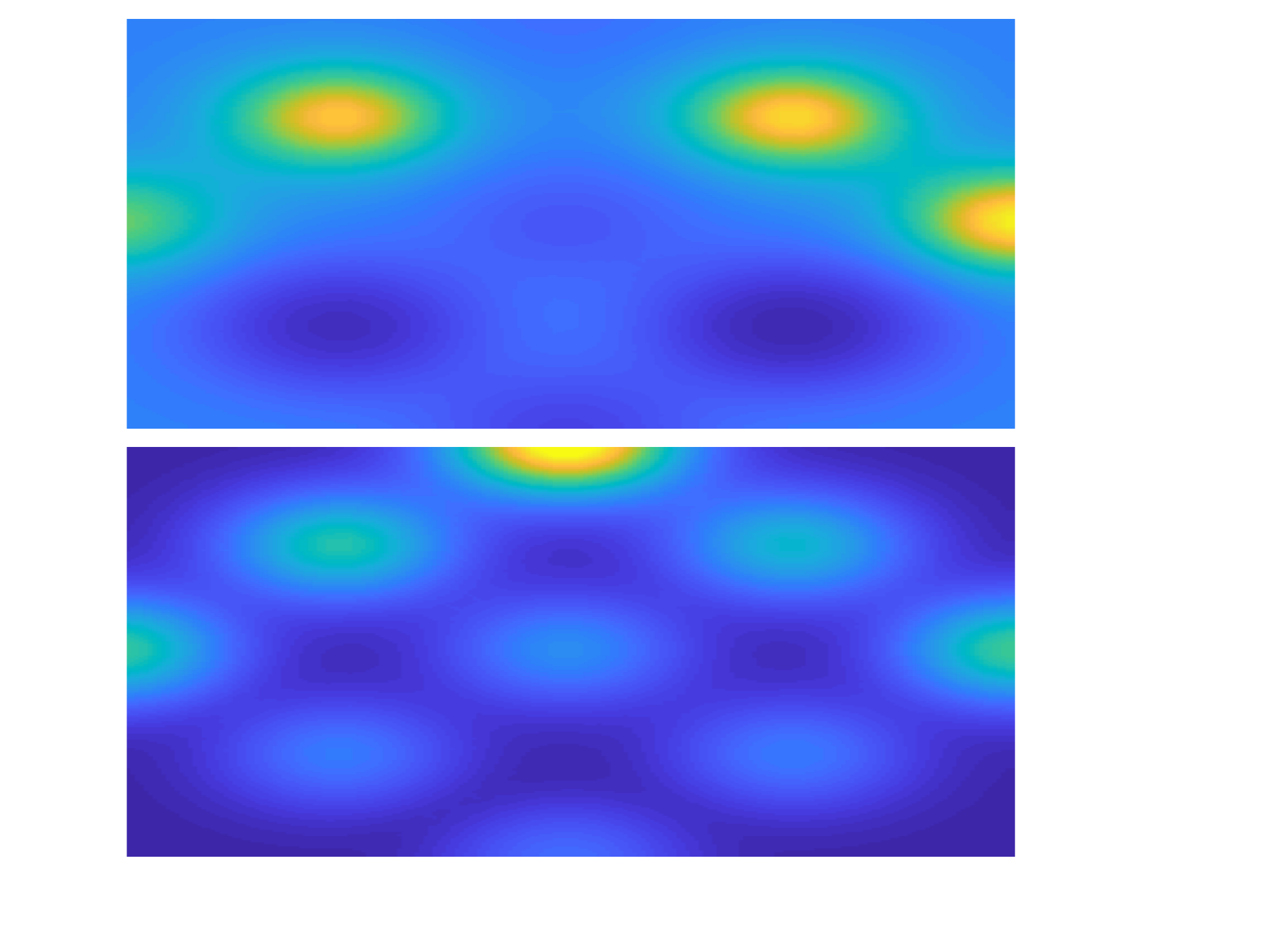}
         \put (15,-3) {\scriptsize NN-RTO, N=50}
  \end{overpic}
  \begin{overpic}[height=6.4cm,width=4.25cm,trim=35 10 45 10, clip=true,tics=10]{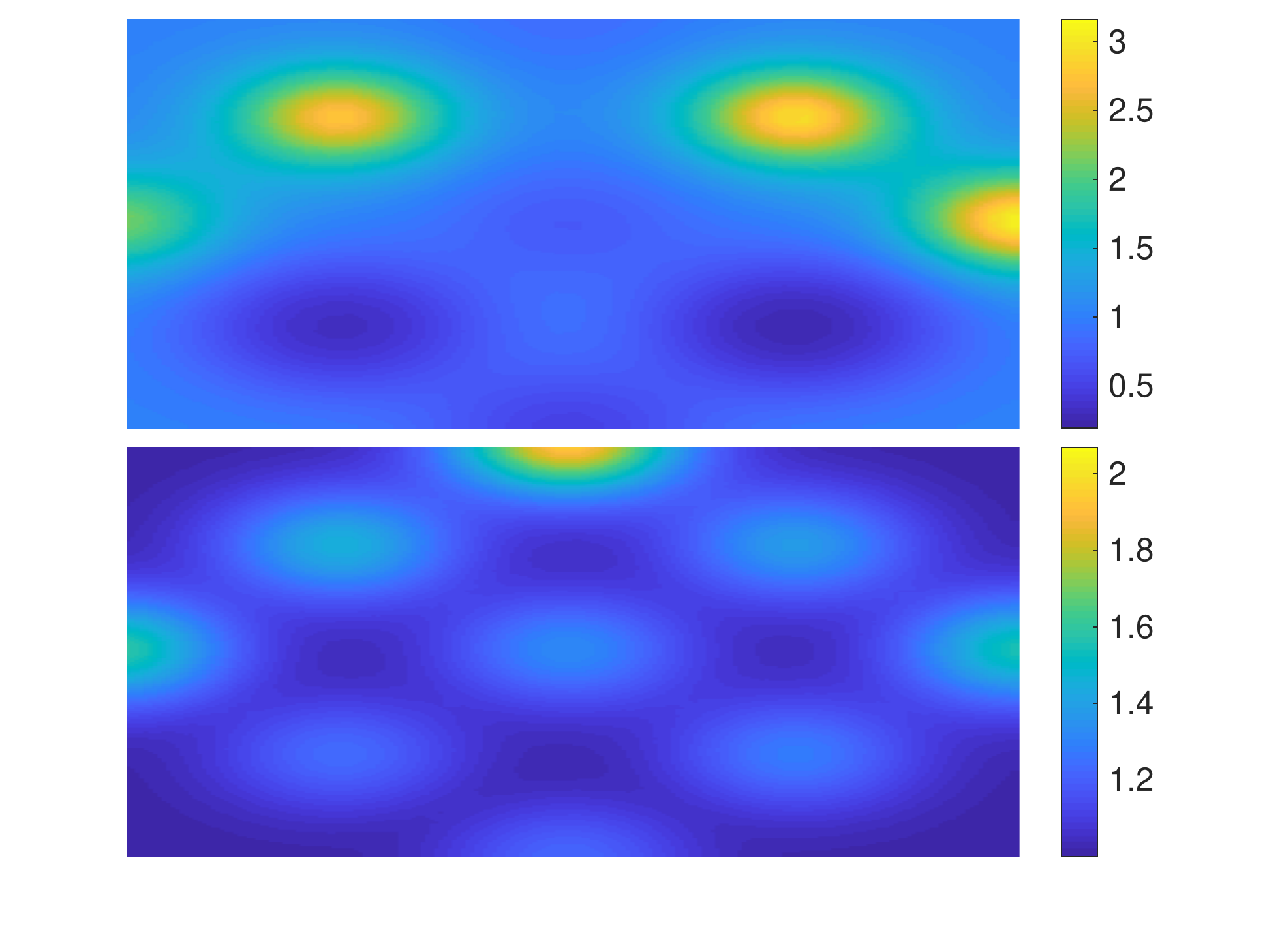}
       \put (15,-3) {\scriptsize NN-RTO, N=100}
  \end{overpic}
  \end{center}
\caption{ Posterior mean (top) and  posterior standard deviation (bottom) of $p(x)$ arising from direct RTO, NN-RTO approach $(N=50, 100)$, respectively. }\label{DNN_sol_eg1}
  \end{figure}

\begin{figure}
\begin{center}
    \begin{overpic}[width=0.7\textwidth,trim= 10 0 40 15, clip=true,tics=10]{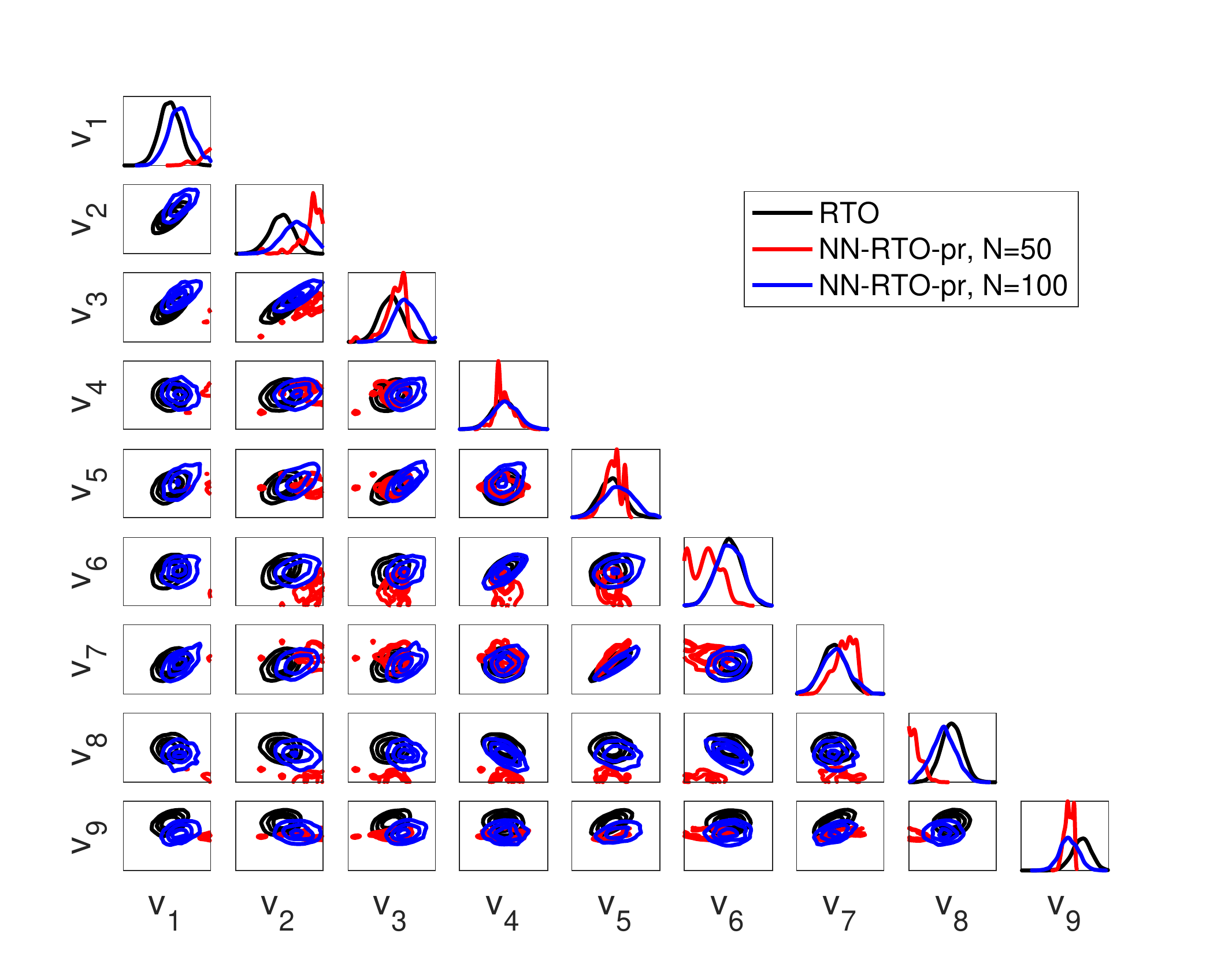}
      \end{overpic}
\end{center}
\caption{One- and two-dimensional posterior marginals of the two parameters for Example 1 using NN-RTO-pr approach. }\label{pr-pos_contour_eg1}

\end{figure}

 \begin{figure}
\begin{center}
    \begin{overpic}[height=6.4cm,width=4.25cm,trim= 35 10 45 5, clip=true,tics=10]{figures/resD_setup0225-eps-converted-to.pdf}
         \put (23,-3) {\scriptsize RTO}
  \end{overpic}
    \begin{overpic}[height=6.4cm,width=4.25cm,trim= 35 10 45 5, clip=true,tics=10]{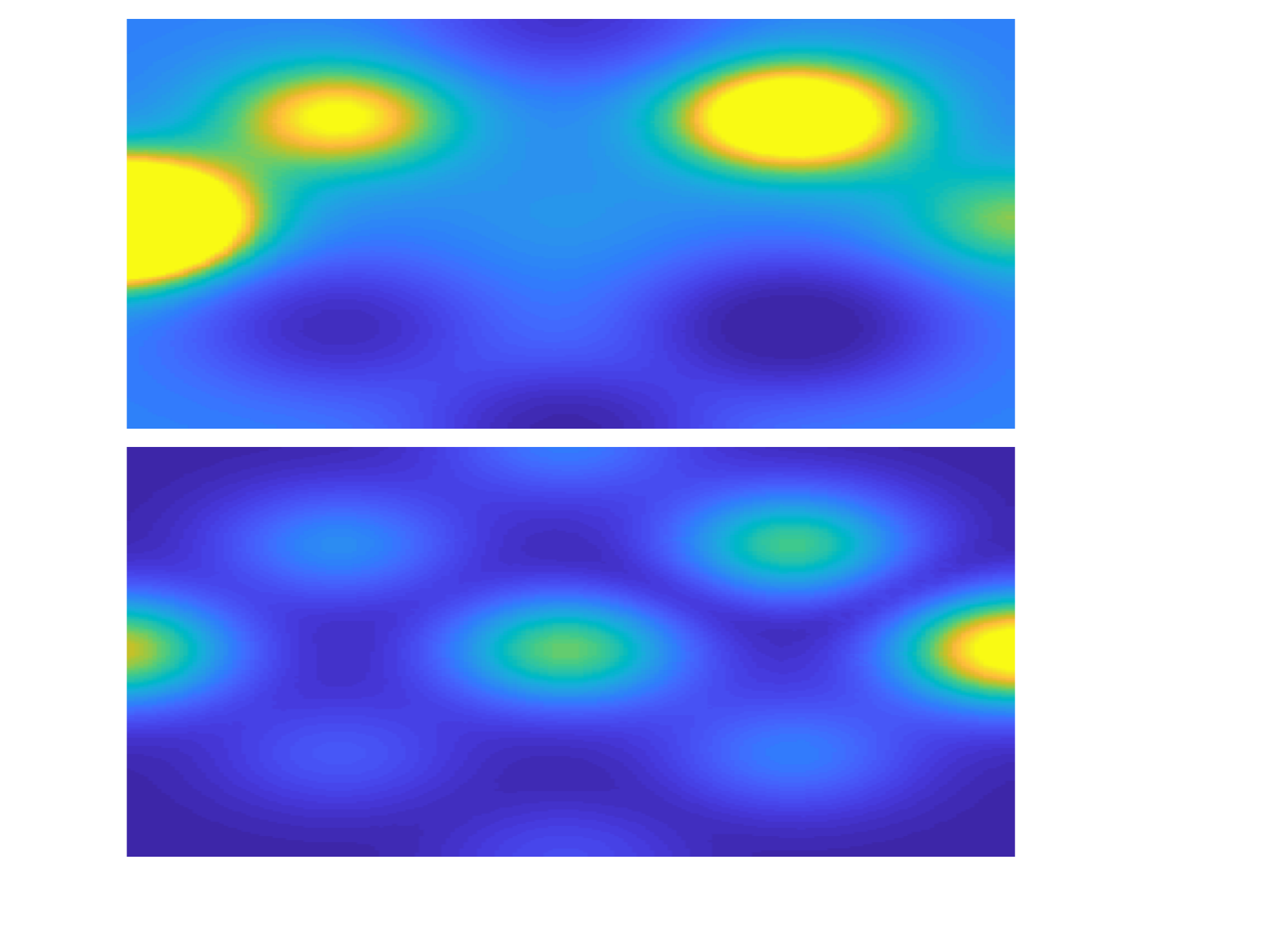}
         \put (15,-3) {\scriptsize NN-RTO-pr, N=50}
  \end{overpic}
  \begin{overpic}[height=6.4cm,width=4.25cm,trim=35 10 45 5, clip=true,tics=10]{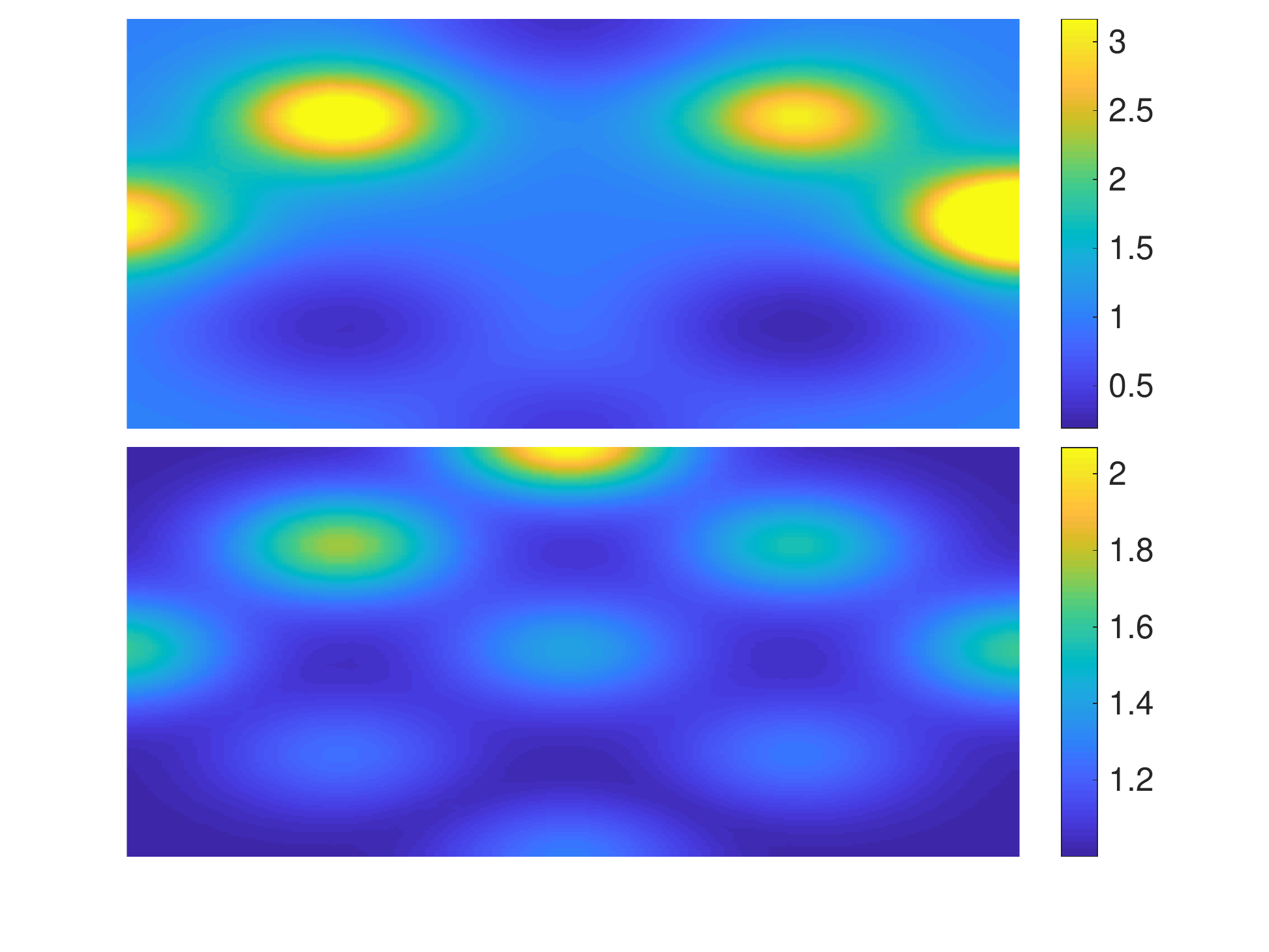}
       \put (15,-3) {\scriptsize NN-RTO-pr, N=100}
  \end{overpic}
  \end{center}
\caption{ Posterior mean (top) and  posterior standard deviation (bottom) of $p(x)$ arising from direct RTO, NN-RTO-pr approach $(N=50, 100)$ , respectively. }\label{pr-DNN_sol_eg1}
  \end{figure}

We first investigate the sampling accuracy of the DNN-RTO algorithm. To this end, we shall construct the DNN surrogate using $N=\{50,100\}$ training points with 3 hidden layers and 40 neurons per layer. Notice that the training points are chosen with the algorithm in Section 4.2. In Fig. \ref{pos_contour_eg1}, we present the marginal distributions of each component of the parameters, and the contours of the marginal distributions of each pair of components. The black lines represent the results generated by the direct RTO approach based on the true forward model evaluations (the reference solution), the red and blue lines represent results of the DNN-RTO with $N=50$  and $100$, respectively. It is clearly seen in Fig. \ref{pos_contour_eg1} that the DNN-RTO algorithm results in a good approximation to the reference solution.  The posterior mean and posterior standard deviation obtained by the DNN-RTO approach and the direct RTO approach are presented in Fig. \ref{DNN_sol_eg1}.  Good agreements between the two algorithms are observed.

Next, we compare the accuracy of our proposed algorithm to that of DNN-RTO-pr (with training points generated according to the prior distribution). Again, we use $N=\{50, 100\}$ training points that are generated by the prior information (rather than the algorithm in Section 4.2) to train the DNN.
The corresponding results are reported in Figs. \ref{pr-pos_contour_eg1} and \ref{pr-DNN_sol_eg1}.  It is shown in the figures that the results using the prior-based DNN approach admits a larger error. By comparing Figs. \ref{pos_contour_eg1} and \ref{DNN_sol_eg1}, we can conclude that the approximation results using DNN-RTO are much more accurate than that of the prior-based "DNN-RTO-pr" approach.

In Table \ref{eg1_time}, we summarize the acceptance probability, the CPU time, the ESS  (min., med., max.), the time-normalized minESS and the speedup factor, for comparing the direct RTO and the DNN-RTO approach. It can be seen that RTO makes 5000 draws in 9806 seconds and obtains an acceptance probability around 0.93. While for DNN-RTO with $N=50$, only 156 seconds are requied to get an acceptance probability around 0.82. DNN- RTO yields 15.4 effective samples per second while RTO is only 0.41 effective draws per second. It is noticed that even when with $N=50$ training points, the DNN-RTO approach can speed up the number of effectively independent samples generated per second by a factor of 37. One can also learn that NN-RTO is much more efficient than DNN-RTO-pr.

% Table==============
 \begin{table}[tp]
      \caption{Comparison of the computational efficiency in Example 1. AP is the acceptance probability, s is the CPU times (second), ESS (min., med., max.), minESS/s is the time-normalized ESS and spdup is the speed up of sampling efficiency measured by minESS/s with RTO as the baseline.  MCMC chain length is 5000 steps.}\label{eg1_time}
  \centering
  \fontsize{6}{12}\selectfont
  \begin{threeparttable}
    \begin{tabular}{ lccccc}
  \toprule
Method                                & AP       &s       & ESS              &  minESS/s &spdup\cr
  \midrule
   RTO                                &  0.93   & 9806   & (4030, 4342, 4460)        &  0.41               & 1      \cr
   DNN-RTO, $N=50$            &  0.82  & 156     & (2409, 2719, 3428)         &  15.4              &37.6  \cr
   DNN-RTO, $N=100$            & 0.80   &  170      & (2612, 2895, 3606)     &  15.4              & 37.6   \cr
   DNN-RTO-pr, $N=50$       &  0.08   & 85      & (3, 5, 77)                      &  0.04              &0.1   \cr
   DNN-RTO-pr, $N=100$        & 0.54   &  169      & (529, 958, 1245)         &  3.1                & 7.6  \cr
  \bottomrule
   \end{tabular}
  \end{threeparttable}
\end{table}
%==============

\subsubsection{Sensitivity to the neural network architecture}
We now investigate the sensitivity of our algorithm with respect to the architecture of the neural networks.  We consider the  error in estimating mean $REM$ and covariance $REC$ of parameters as
\begin{eqnarray*}
REM= \frac{\|\bar{\kappa}-\kappa^{\dag}\|_{\infty}}{\|\kappa^{\dag}\|_{\infty}},
\end{eqnarray*}
and
\begin{eqnarray*}
REC= \frac{\|cov(\kappa)-cov(\kappa^{\dag})\|_{F}}{\|cov(\kappa^{\dag})\|_{F}},
\end{eqnarray*}
where $\kappa^{\dag}$ and $cov(\kappa^{\dag})$ are the ``true" posterior mean and covariance arising from direct RTO,  and $\|\cdot\|_F$ denotes the Frobenius norm.

To verify the sensitivity of the proposed method with respect to the structures of the $\NN$,  we test several constant values choosing with different depth $L\in\{1,2,3,4\}$ and width $d_k\in\{20,40,60,80\}$. With these settings, we  run Algorithm \ref{alg:DNN-RTO-MH} using the $N=100$ training points.  The corresponding numerical results  are presented in Table \ref{rel_L}, in which the corresponding CPU times are also presented.  As shown in this table, the computational results for DNN-RTO with different depth $L$ and width $d_k$ are almost the same. Table \ref{rel_N} shows the  estimate errors $REM$, $REC$ and the CPU times $s$ with respect to the number of training points obtained with various values of $d_k$ and $L=3$ hidden layers. Overall, we observe that the numerical results remain robust for all neural network architectures considered.

\begin{table}
  \centering
 % \fontsize{8}{10}\selectfont
\begin{tabular}{|c|c|c|c|c|}
\hline
\diagbox{L}{$d_k$} & 20& 40 & 60 & 80 \\
\hline
1 & (0.019,0.102,94) & (0.017,0.046,121) & (0.022,0.052,128) & (0.017,0.033,137) \\
\hline
2 & (0.026,0.510,150) & (0.022,0.095,166) & (0.010,0.057,195) & (0.019,0.038,202)\\
\hline
3 & (0.029,0.949,139) & (0.029,0.068,210) & (0.018,0.062,246) & (0.016,0.077,270) \\
\hline
4 & (0.402,0.918,175) & (0.029,0.074,241) & (0.055,0.100,286) & (0.049,0.066,293)\\
\hline
\end{tabular}
  \caption{Example 1:  The estimate errors $REM$, $REC$ and the CPU times $s$ obtained using  NN-RTO approach with various values of $L$, $d_k$ and $N=100$ training points.}\label{rel_L}
    \centering
 % \fontsize{8}{10}\selectfont
\begin{tabular}{|c|c|c|c|c|}
\hline
\diagbox{N}{$d_k$} & 20& 40 & 60 & 80 \\
\hline
50 &  (0.069,1.150,155) & (0.058,0.097,191) & (0.019,0.072,228) & (0.022,0.033,261) \\
\hline
100 & (0.029,0.949,139) & (0.029,0.068,210) & (0.018,0.062,246) & (0.016,0.077,270)\\
\hline
150 & (0.018,0.136,193) & (0.016,0.034,209) & (0.015,0.054,227) & (0.019,0.038,255) \\
\hline
200 &  (0.018,0.134,117) & (0.015,0.032,210) & (0.019,0.063,225) & (0.018,0.050,256)\\
\hline
\end{tabular}
  \caption{The estimate errors $REM$,$REC$ and the CPU times $s$ obtained using  NN-RTO approach with various values of $N, d_k$ and $L=3$ hidden layers. }\label{rel_N}
\end{table}

\subsection{Example 2: a high dimensional inverse problem}

% Table==============
 \begin{table}[tp]
      \caption{Comparison of the computational efficiency in Example 2. AP is the acceptance probability, s is the CPU times (second), ESS(min., med., max.),  minESS/s is the time-normalized ESS and spdup is the speed up of sampling efficiency measured by minESS/s with RTO as the baseline.  MCMC chain is 5000 steps.}\label{eg2_time}
  \centering
  \fontsize{6}{12}\selectfont
  \begin{threeparttable}
    \begin{tabular}{ lccccc}
  \toprule
Method                                & AP       &s       & ESS              &  minESS/s &spdup\cr
  \midrule
   RTO                                       &  0.54   & 2271   & (1005, 1249, 1668)          &  0.44              & 1      \cr
   DNN-RTO, $N=1000$            &  0.35   & 96    & (1050, 1983, 2259)       &  10.9                 &25    \cr
   DNN-RTO, $N=1500$            & 0.35   & 88      & (1036, 1264, 1470)        &  11.8               & 27   \cr
  \bottomrule
   \end{tabular}
  \end{threeparttable}
\end{table}
%==============

\begin{figure}
\begin{center}
    \begin{overpic}[width=0.45\textwidth,trim= 30 10 10 25, clip=true,tics=10]{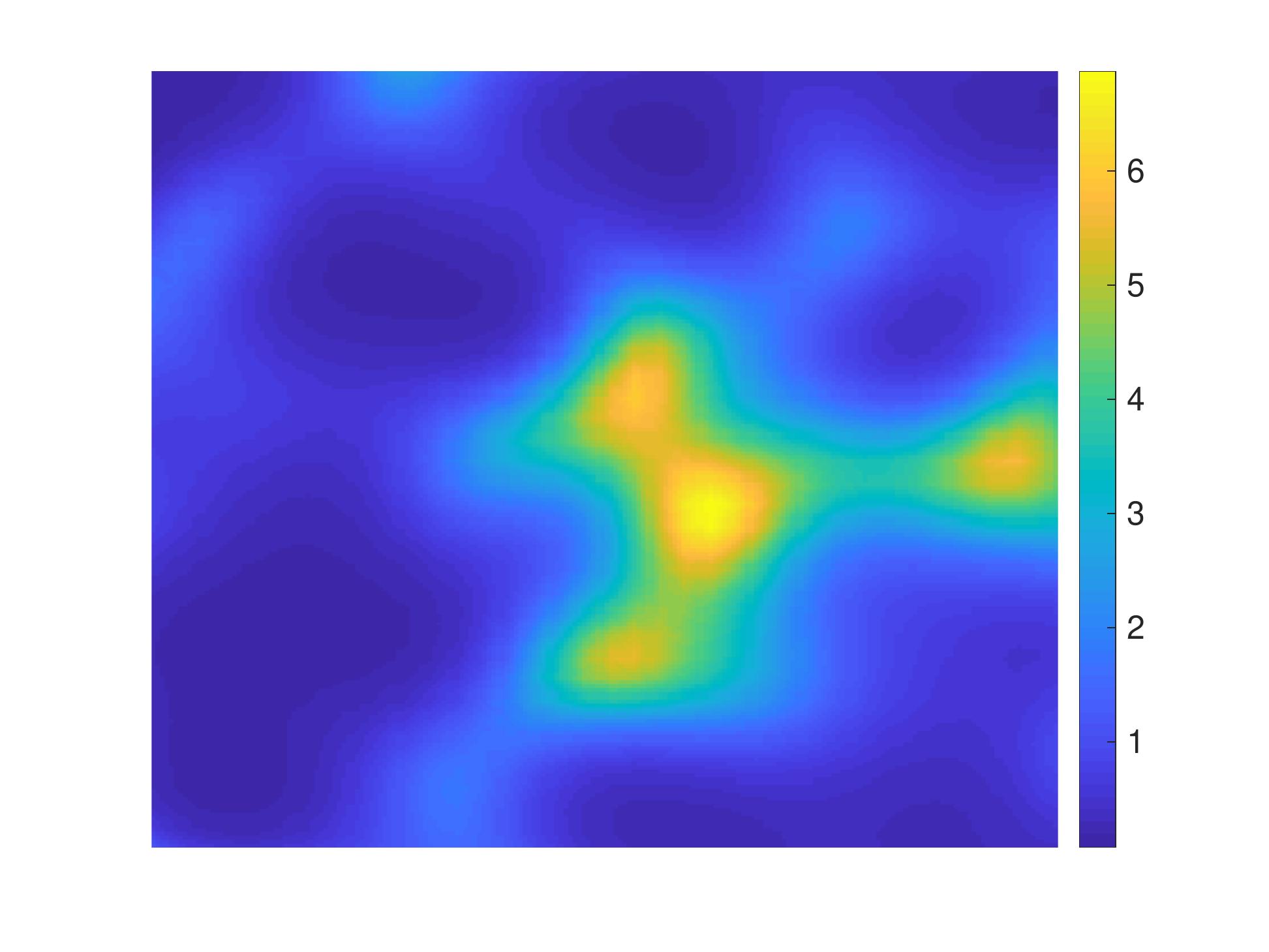}
      \end{overpic}
          \begin{overpic}[width=0.45\textwidth,trim= 30 10 10 25, clip=true,tics=10]{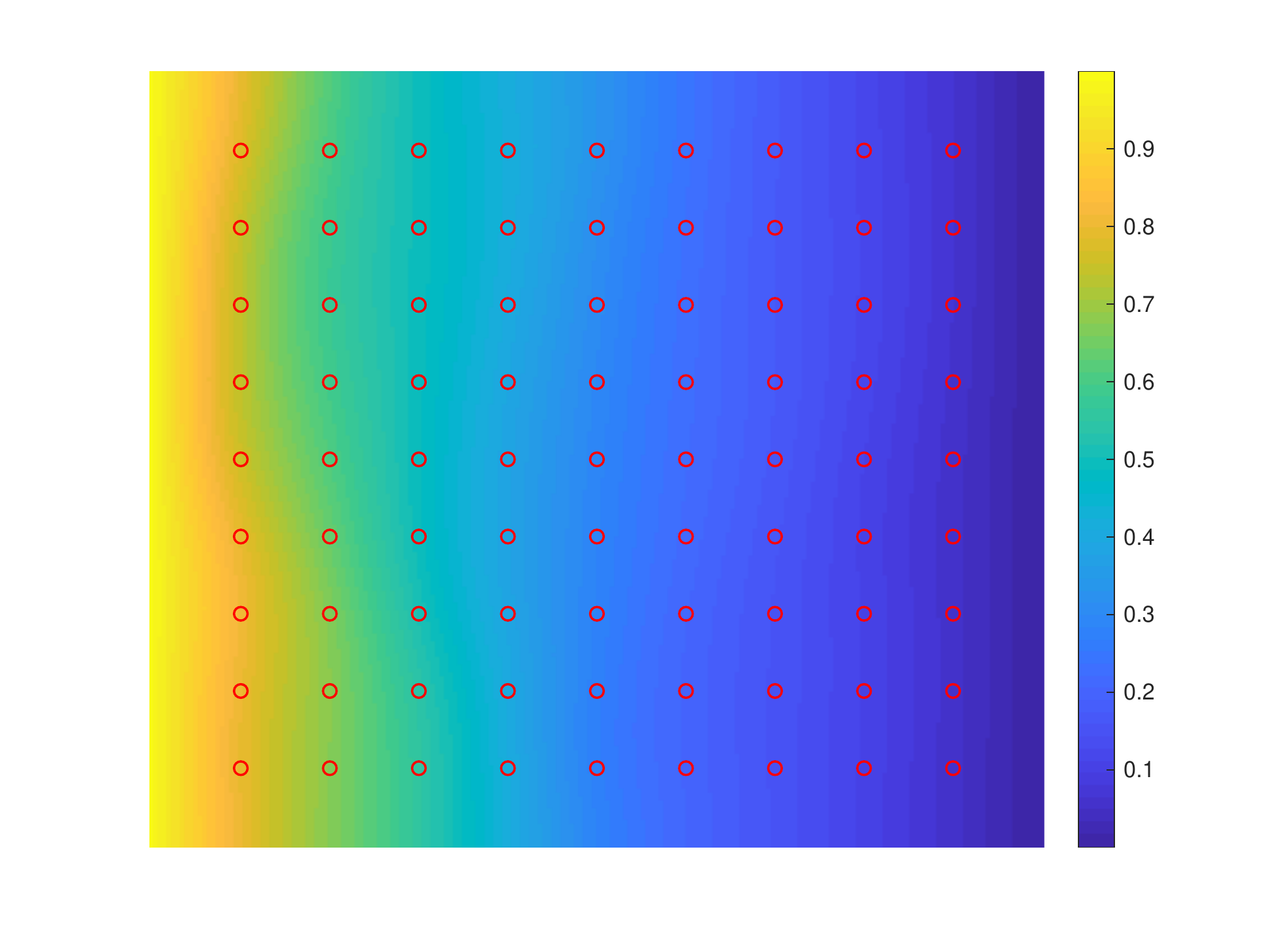}
      \end{overpic}
\end{center}
\caption{The true permeability $\kappa(x)$(left) and the model outputs (right) for elliptic PDE inverse problems in example 2. }\label{exact_eg2}
\end{figure}

In the second example, we consider a log-diffusivity field $\log\kappa(x)$ that is endowed with a Gaussian process prior, with zero mean and an isotropic kernel:
\begin{equation*}
C(x_1,x_2)=\sigma^2 \exp\Big(-\frac{\|x_1-x_2\|^2}{2l^2}\Big).
\end{equation*}
Here we set variance $\sigma^2=1$ and $l = 0.1$. This prior allows the field to be easily parameterized with a Karhunen-Loeve expansion:
\begin{equation}
p(x; v) \approx \sum^{n}_{i=1} v^i \sqrt{\lambda_i} \phi_i(x),
\end{equation}
where $\lambda_i$ and $\phi_i(x)$ are the eigenvalues and eigenfunctions, respectively, of the integral operator on $[0,1]^2$ defined by the kernel $C$, and the parameter $v=(v^1,\cdots, v^n)$ are endowed with independent standard normal priors, $v^i \sim N(0,1)$. These parameters then become the targets of inference. In particular, we truncate the Karhunen-Loeve expansion with $n=120$ modes.  The true solution  $\kappa(x)$  used to generate the test data are shown in Fig.\ref{exact_eg2}.  The measurement sensors of $p$ are evenly distributed over $[0.1,0.9]^2$ with grid spacing 0.1, i.e., $d \in \R^{81}$.  The observational errors are taken to be additive and Gaussian:
\begin{equation*}
d_j = p(x_j;v) +\xi_j, \quad j=1,\cdots,81,
\end{equation*}
with $\xi_j \sim N(0,0.05^2)$.  In this example, four hidden layers and 80 neurons per layer are used in $\NN$.

\begin{figure}
\begin{center}
    \begin{overpic}[height=6.4cm,width=4.25cm,trim= 35 10 45 5, clip=true,tics=10]{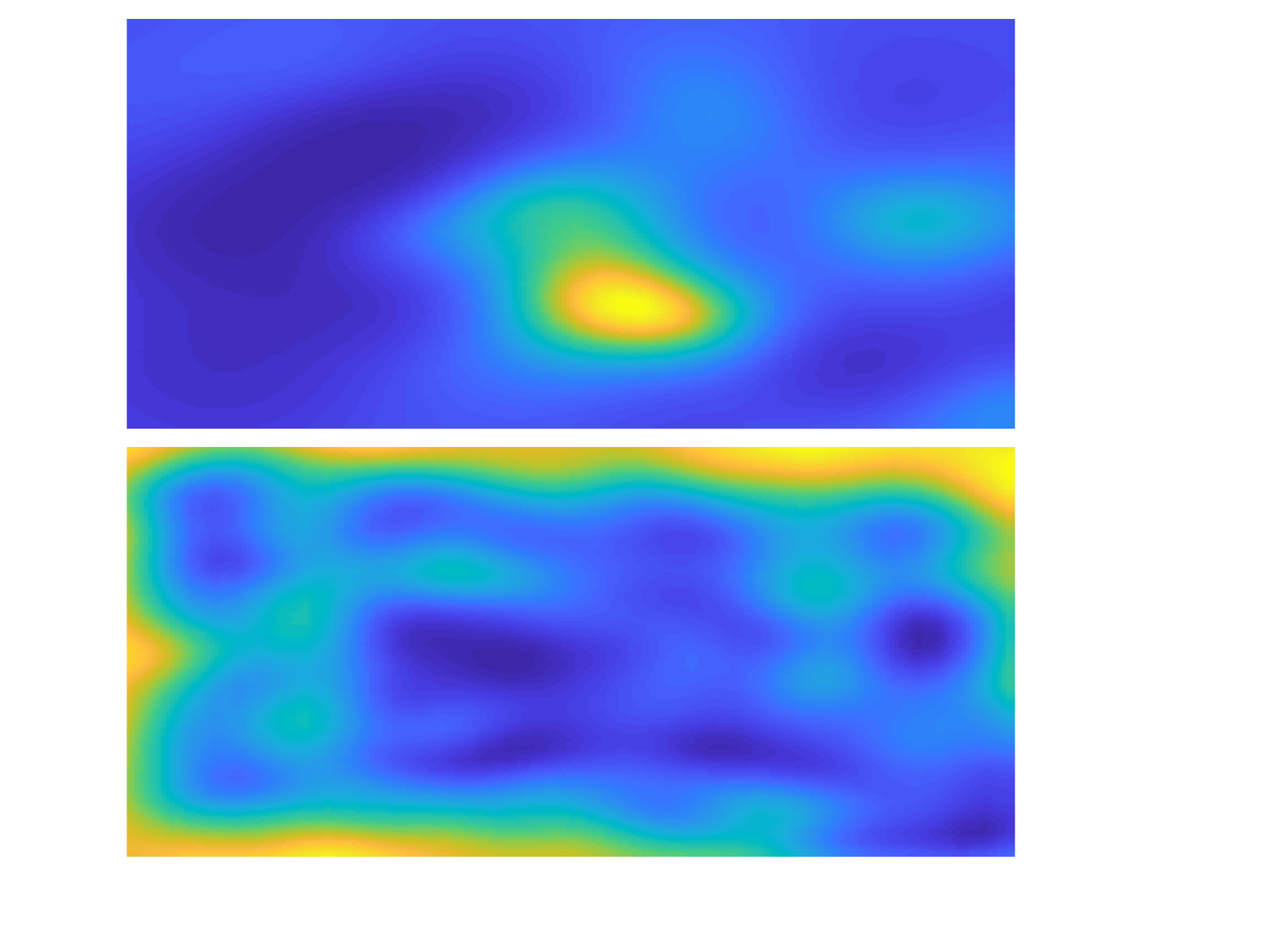}
         \put (23,-3) {\scriptsize RTO}
  \end{overpic}
    \begin{overpic}[height=6.4cm,width=4.25cm,trim= 35 10 45 5, clip=true,tics=10]{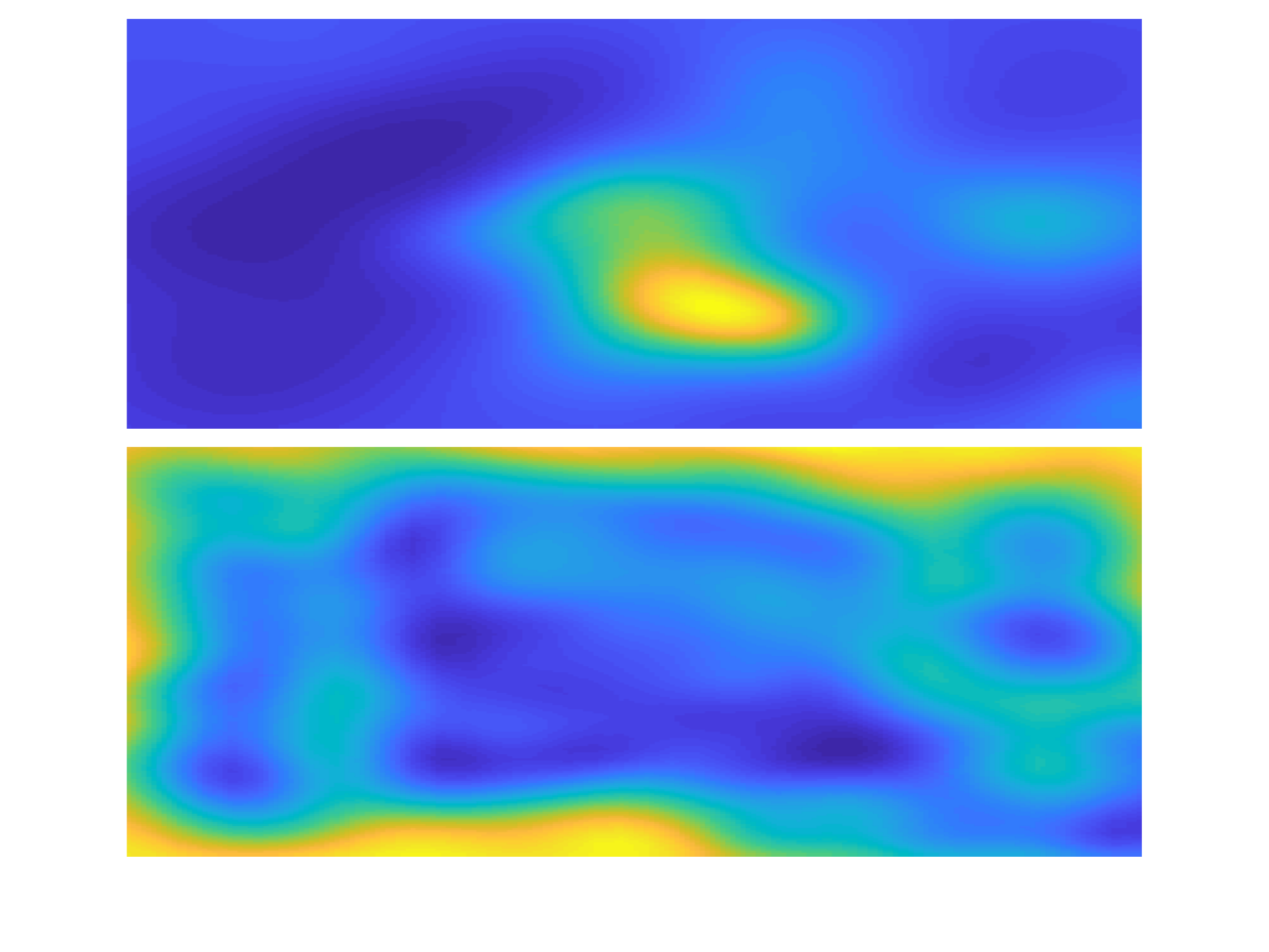}
         \put (15,-3) {\scriptsize NN-RTO, N=1000}
  \end{overpic}
  \begin{overpic}[height=6.4cm,width=4.25cm,trim=35 10 45 5, clip=true,tics=10]{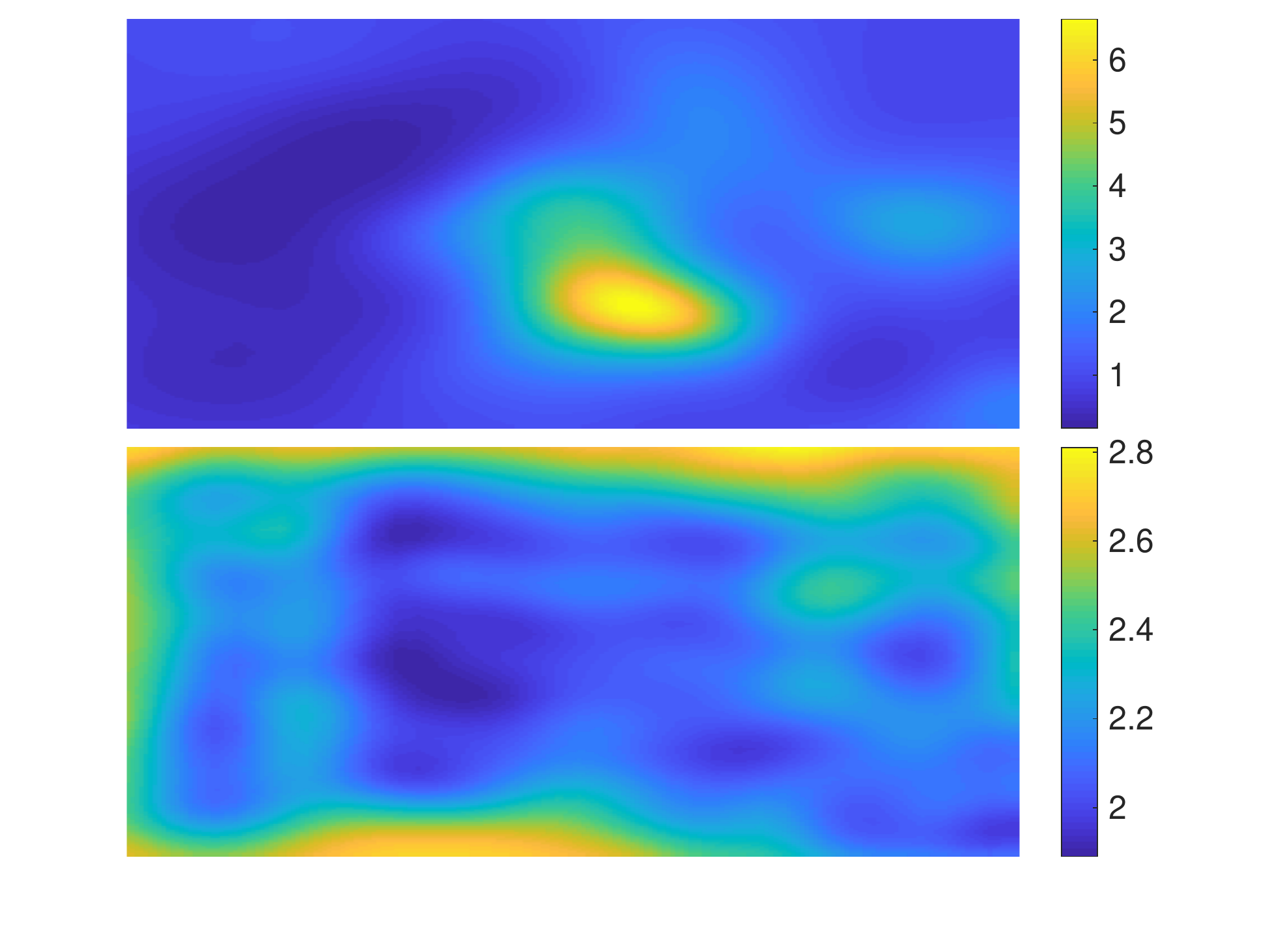}
       \put (15,-3) {\scriptsize NN-RTO, N=1500}
  \end{overpic}
  \end{center}
\caption{ Posterior mean (top) and  posterior standard deviation (bottom)  arising from direct RTO, NN-RTO approach ($N=1000$) and prior-based NN-RTO ($N=1500$), respectively. } \label{pmean_eg2}
  \end{figure}

We ran  DNN-RTO  with $N=\{1000, 1500\}$ training points and generated an MCMC chain of length 5000.  As shown in Fig. \ref{pmean_eg2}, the posterior mean and posterior standard deviation estimated by DNN-RTO are closer to the reference solution (that is obtained by direct RTO approach). We compare the sampling efficiency of different algorithms, and the result is summarized in Table \ref{eg2_time}. Again, we observe a similar result in the raw ESS when comparing DNN-RTO algorithms with the direct RTO, but an increase in efficiency due to the computational time cut by DNN surrogate.  The computation used 96 seconds to produce a minESS of  about $10^3$ for DNN-RTO with $N=1000$.  Notice that more than an order of magnitude of improvement is observed for DNN-RTO compared to RTO.  All the above discussion confirm that DNN-RTO is advantageous in sampling efficiency.

\section{Summary} \label{sec:summary}

In this paper, we present a new strategy, namely the DNN-RTO algorithm, to accelerate the original RTO-MH algorithm.  One of the key components of our DNN-RTO algorithm is a goal-oriented strategy for choosing the training points from a local Gaussian measurement.  Since the numerical accuracy of the Bayesian inverse problems is mainly concerned in a posterior density region, our DNN-surrogate requires very few training points to achieve the same level of accuracy compared with a prior-based DNN-surrogate. To demonstrate the accuracy and efficiency of the proposed algorithm, a benchmark example to infer the permeability field for elliptic PDEs with synthetic data is tested. The numerical results show that the DNN-RTO is able to accelerate RTO sampling by up to several orders of magnitude. We believe the approach in this work will be promising in dealing with high dimensional BIPs and/or BIPs with limited regularity.

\end{document}